\documentclass{CSMA2026}
\title{Efficient multidisciplinary design via Bayesian optimization}
\author{N. Bartoli$^{1,2}$, T. Lefebvre$^{1,2}$, R. Lafage$^{1,2}$, P. Saves$^{3}$, Y. Diouane$^{4}$, J. Morlier$^{5}$}
\address{%
$^1$ DTIS, ONERA, Université de Toulouse, 31000, France\\
$^2$ Fédération ENAC ISAE-SUPAERO ONERA, Université de Toulouse, 31000, Toulouse, France \\
$^3$IRIT, UMR 5505 CNRS, Université Toulouse Capitole, Université de Toulouse, Toulouse, France\\
$^4$ GERAD and Department of Mathematics and Industrial Engineering, Polytechnique Montr\'eal, Montréal, QC, Canada \\
$^5$ ICA, Université de Toulouse, ISAE–SUPAERO, INSA, CNRS, MINES
ALBI, UPS, Toulouse, France
}

\usepackage{natbib} 

\begin{document}
\maketitle

\begin{abstract}
This study introduces SEGOMOE, a Bayesian optimization tool for optimizing complex, computationally expensive systems, especially in aeronautics. It efficiently handles mixed design variables (continuous, discrete, categorical, hierarchical) using adaptive Gaussian process models. SEGOMOE combines expert models to address nonlinearities in objectives and constraints, leveraging the open-source Surrogate Modeling Toolbox (SMT).
The tool supports multi-fidelity data and solves both single- and multi-objective problems, including hidden constraints and high-dimensional decomposition. Validated through benchmarks and real-world aeronautical applications, SEGOMOE proves to be robust and versatile for tackling multidisciplinary challenges.
\end{abstract}

\section{Context}
Modern engineering challenges, especially in aerospace and aeronautics, necessitate dedicated optimization processes to achieve optimal performance, often in contexts where new architectural concepts are explored to meet ambitious goals such as sustainability. Multidisciplinary Design Optimization (MDO) aims to find the best design by considering the trade-offs and dependencies between various disciplines (e.g., aerodynamics, structure, propulsion)~\cite{Martins2021Engineering}.
The nature of these problems poses significant challenges for traditional optimization methods. MDO relies heavily on physics-based simulations to accurately predict the performance of design candidates, particularly when historical data or expertise is lacking for novel architectures. Consequently, the objective functions and constraints are treated as expensive-to-evaluate black-box functions. Calculating derivatives is often impractical due to the complexity or high computational costs associated with methods like finite differencing~\cite{audet2017derivative}.
Bayesian Optimization~(BO) offers a particularly effective solution within this context~\cite{ Jones1998Efficient,garnett2023bayesian}. As a derivative-free global optimization approach, BO is designed to minimize the total number of expensive black-box function calls required to locate the optimum in MDO engineering domains~\cite{Bartoli2019AdaptiveModeling, Tfaily2024BayesianHC} and System Architecture Optimization~(SAO)~\cite{Bussemaker2021JetEngine,  Bussemaker2024SAOStrategies}. BO overcomes the cost difficulty by employing probabilistic surrogate models, typically Gaussian Processes (GP)~\cite{Rasmussen2006Gaussian, Kennedy2002Bayesian}, which provide both a mean prediction and uncertainty quantification. The iterative BO framework, such as SEGOMOE~\cite{Bartoli2019AdaptiveModeling, Priem2020EfficientApplication} for Super Efficient Global Optimization with Mixtures of Experts, adaptively guides the search by maximizing an inexpensive acquisition function, ensuring an efficient balance between exploration and exploitation~\cite{Jones1998Efficient, garnett2023bayesian}. This paper reviews the adaptation of BO techniques to address the multifaceted complexities of real-world MDO problems, including constrained multi-objective formulations~\cite{Charayron2023MultiFidelityMO, Grapin2022ConstrainedMOBO, donelli2025concurrent}, the management of non-linear and hidden constraints~\cite{Tfaily2024BayesianHC, Priem2020UTB, Muller2019Surrogate}, and the modeling of mixed and hierarchical design variables~\cite{Bussemaker2024SAOStrategies,Saves2024SMT2, Saves2023MixedCategorical}, with applications drawn from aeronautical engineering such as the simultaneous trajectory and design optimization of drones~\cite{Fernandez2023Assessment} and turbomachinery optimization~\cite{Bussemaker2021JetEngine, Pretsch2025CooperativeComponentsBO}.

\section{BO for multidisciplinary and black-box design optimization}
The design of complex systems, such as aircraft architectures, often requires MDO, seeking optimal solutions by managing trade-offs among coupled disciplines~\cite{Martins2021Engineering,Lambe2012Extensions, Cramer1994ProblemFormulation}. MDO involves significant computational effort, as it relies heavily on costly physics-based numerical simulations, making the objective and constraint functions expensive-to-evaluate black-box functions~\cite{Jones1998Efficient, Kennedy2002Bayesian,LeGratiet2013MultiFidelity}. Furthermore, derivatives are typically unavailable or impractical to compute when mixed variables are involved~\cite{Bussemaker2024SAOStrategies,Saves2024HighDimMixed}.

\subsection{The BO framework}
BO is a derivative-free global optimization strategy designed to minimize the total number of expensive black-box function calls required to find the optimum~\cite{Jones1998Efficient, garnett2023bayesian,Frazier2018Tutorial}.
BO employs a surrogate model, predominantly the GP (or Kriging model)~\cite{Rasmussen2006Gaussian, Kennedy2002Bayesian}, which models the unknown function by providing both a mean prediction ($\mu$) and an associated uncertainty (variance $\sigma^2$).
The process is iterative~\cite{Jones1998Efficient}:
\begin{enumerate}
    \item \textbf{Model training:} The GP model is built or updated using existing data, denoted as the Design of Experiments (DoE)~\cite{Rasmussen2006Gaussian, Sacks}.
    \item \textbf{Acquisition function optimization:} An inexpensive acquisition function (e.g., Expected Improvement, Upper Confidence Bound, Watson-Barnes 2, WB2s) is maximized to determine the next point to evaluate~\cite{Jones1998Efficient, Bartoli2019AdaptiveModeling,Sasena02flexibility}. This function balances exploitation (choosing promising points near the current minimum) and exploration (investigating areas of high uncertainty).
    \item \textbf{Evaluation and update:} The chosen point is evaluated using the costly black-box simulation, and the DOE is augmented for the next iteration.
\end{enumerate}
The GP models are implemented in the open-source Surrogate Modeling Toolbox (SMT)\footnote{https://github.com/SMTorg/smt}~\cite{Saves2024SMT2, Bouhlel2019PythonSMT}. 

\subsection{BO extensions of complex optimization problems}
Different extensions have been added in the initial process to be applied to more complex problems involving multiple objectives, non linear and hidden constraints, multiple fidelity models with different accuracy and CPU time, and/or mixed variables. The following paragraphs describe the proposed contributions.
\paragraph{Mono- and multi-objective formulation}
MDO often involves optimizing multiple conflicting criteria, necessitating Multi-Objective Bayesian Optimization (MOBO) to find the Pareto Front (PF)~\cite{Deb2002NSGAII, Emmerich2006MultiObjective}.  Extensions of the SEGOMOE algorithm for multi-objective problems have shown that they can obtain Pareto fronts with significantly fewer evaluations compared to methods like NSGA-II~\cite{Grapin2022ConstrainedMOBO, Deb2002NSGAII}. 
The main modifications were made to the acquisition functions, using hypervolume improvement with various strategies~\cite{Grapin2022ConstrainedMOBO}, and an additional post-processing step was added to construct the predicted Pareto front (based on the final surrogate model). This step aims to both complete and densify the final Pareto front and compare it with the recomputed Pareto front (derived from the black-box simulations) to assess the convergence of the process.
\paragraph{Mono- and multi-fidelity formulation}
In a MDO context, evaluating objectives and constraints often involves complex, computationally expensive black-box models.
Multi-fidelity approaches can help reduce these costs by combining models with varying levels of accuracy and computational expense.
Specialized MOBO frameworks, such as MFMO-SEGO for Multi-Fidelity \& Multi-Objective Super Efficient Global Optimization, extend the methodology to use multiple fidelity levels~\cite{Charayron2023MultiFidelityMO,LeGratiet2013MultiFidelity} in a two-step approach, where first the next interesting point to be evaluated is found, and then the fidelity level is determined~\cite{MFEGO_AIAA:19,cordelier2025}.

\paragraph{Constrained optimization and feasibility}
Engineering problems involve non-linear, black-box constraints~\cite{Bartoli2019AdaptiveModeling,le2024taxonomy}. The SEGOMOE framework models both objectives and constraints (inequality and equality) using GP surrogates~\cite{Bartoli2019AdaptiveModeling, Priem2020UTB}.
To handle complex constraints, the Upper Trust Bound (UTB) feasibility criterion was developed~\cite{Priem2020UTB}, dynamically combining the GP mean prediction and its associated uncertainty during the optimization process.
This approach encourages exploration near constraint boundaries or in uncertain regions, which is crucial for maximizing the objective function under strict constraints.

\paragraph{Handling hidden constraints}
Hidden constraints are areas where the expensive simulation fails to converge, crashes, or returns an error (NaN)~\cite{Muller2019Surrogate, Forrester2006MissingData}. These failures can arise from solver divergence~\cite{Bussemaker2024SurrogateHC} or physical infeasibility, such as the failure of a landing gear sizing code in certain configurations.
Strategies include predicting the failure region using supervised machine learning (ML) classifiers, such as Random Forest Classifiers (RFC), K-Nearest Neighbors (KNN), or GP classifiers~\cite{Tfaily2024BayesianHC}.
A novel approach is the Feasibility Enhanced Expected Improvement (EFIFE) acquisition function, which integrates output from ML classifiers modeling the hidden constraints~\cite{Tfaily2024BayesianHC}. EFIFE aims to reduce the classifier's negative impact on the exploration part of the criterion, leading to improved objective function results with fewer evaluations in aircraft conceptual design~\cite{Tfaily2024BayesianHC}.

\paragraph{Mixed and hierarchical design variables}
System Architecture Optimization (SAO) involves mixed variables (continuous, integer, and categorical)~\cite{Saves2024SMT2,Saves2023MixedCategorical, Saves2024HighDimMixed}. Hierarchical variables occur when the activeness of one variable depends on the choice of another~\cite{Bussemaker2024SAOStrategies, Saves2024SMT2}.
To model these complex spaces, specialized Gaussian process kernels are implemented in the open-source Surrogate Modeling Toolbox SMT 2.0~\cite{Saves2024SMT2, Bouhlel2019PythonSMT}. SMT 2.0 supports kernels such as Gower Distance (GD), Continuous Relaxation (CR), Homoscedastic Hypersphere (HH), Exponential Homoscedastic Hypersphere (EHH), and Hierarchical (HIER)~\cite{Saves2024SMT2, Saves2023MixedCategorical}. Since the SEGOMOE framework uses GP models from SMT, it greatly facilitates the use of mixed and hierarchical variables~\cite{saves2025modelinghierarchicalspacesreview}.

\paragraph{High-dimensional optimization }
For High-Dimensional Optimization (HDBO) ($d \gg 20$), dimension reduction techniques are integrated~\cite{Saves2024HighDimMixed, Priem2024HighDimBO_AIAAJ}. Kriging with Partial Least Squares (KPLS) extends PLS regression to integrate it into GP models~\cite{Saves2024HighDimMixed,Bouhlel2018EGOHighDim}. KPLS constructs a low-dimensional embedding subspace, which is crucial for handling high-dimensional mixed-categorical MDO by reducing the number of effective dimensions~\cite{Saves2024HighDimMixed, Bouhlel2018EGOHighDim}. Another strategy based on cooperation is also proposed: the overall optimization task is decomposed into lower-dimensional component subproblems. Component interactions are fully taken into account by a sequential cooperative procedure named cooperative components
BO (CC-BO)~\cite{Pretsch2025CooperativeComponentsBO,zhan2024cooperative}.
Some combinations of all these features are illustrated in the following section with numerical applications. 
\section{Applications of Bayesian optimization in aeronautics}




\subsection{DRAGON hybrid electric aircraft concept}
\begin{figure}[htp]
\begin{minipage}[b]{.6\linewidth}
\centering
\hspace{-1.25cm}
{
\includegraphics[height=4.5cm,,width=7cm]
{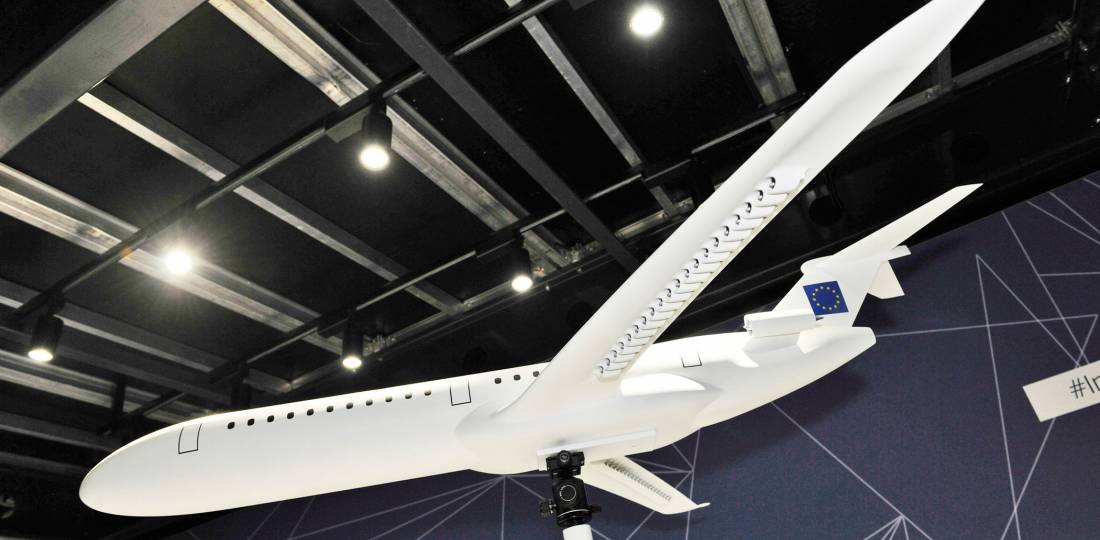}
}
\end{minipage}
\begin{minipage}[b]{.4\linewidth}
\centering 
\hspace{-1.5cm}
{
\includegraphics[height=5cm,width=7cm]{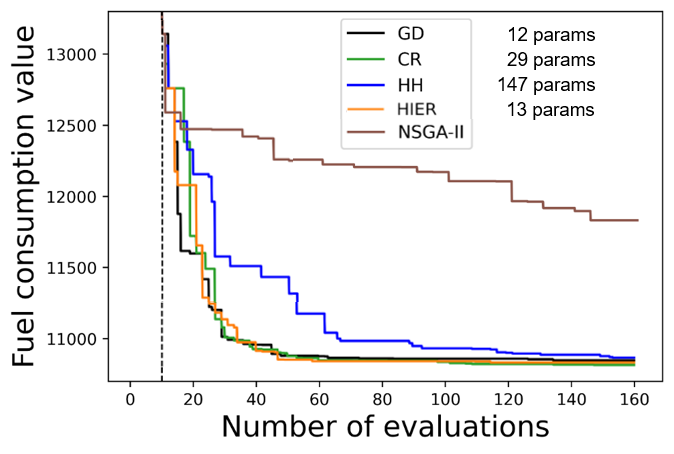}
}
\end{minipage}
\caption{Optimization results for the ``\texttt{DRAGON}'' (mock-up on the left) aircraft~\cite{SciTech_cat}  for 10 DoE of 10 points with the convergence curves: medians of 10 runs (on the right).}
\label{fig:dragon}
\end{figure}

Optimization of the DRAGON concept~\cite{Schmollgruber2019DRAGONExploration} (see Figure~\ref{fig:dragon}) involves six quantities of interest derived from a multidisciplinary analysis involving different disciplines based on the Future Aircraft Sizing Tool with Overall Aircraft Design (FAST-OAD)~\cite{David_2021}. Table~\ref{tab:dragon} lists the objective function to minimize (fuel mass), the five nonlinear constraints, and the different variables describing the geometry (10 continuous design variables) and the propulsion system via some categorical variables. 
\begin{table}[htp]
\small
\centering
\vspace*{-0.3cm}
\caption{Definition of the ``\texttt{DRAGON}'' optimization problem.}
\begin{tabular}{lllrr}
& Function/variable & Nature & \# & Range\\
\hline
\hline
Minimize & Fuel mass & cont & 1 &\\
\hline
w.r.t 
 & \multicolumn{2}{l}{Total  continuous variables} & 10 & \\
 \hline
& \mbox{Architecture} (levels) & cat & 17  & \{1,2, \ldots,16,17\} \\
& \mbox{Turboshaft layout}  (levels) & cat & 2  & \{1,2\} \\
 & \multicolumn{2}{l}{Total categorical variables} & 2 & \\
  \multicolumn{3}{l}{\textbf{Categorical strategy: total number of variables in relaxed dimension }} & {\textbf{29}} & \\
 \hline
 \hline
 & \mbox{Number of generators} (3 levels) & hier & 1  & \{2,4,6\} \\
  & \mbox{Number of motors} (9 levels) & hier & 1  & \{8,12,\ldots,40\} \\
& \mbox{Turboshaft layout}  (levels) & cat & 2  & \{1,2\} \\

 \multicolumn{3}{l}{\textbf{Categorical \& hierarchical strategy: total number of variables}} & {\textbf{14}} & \\
  \hline
subject to & Wing span \textless  \ 36   ($m$)  & cont & 1 \\
 & Take Off Field Length \textless  \ 2200 ($m$) & cont & 1 \\
 & Wing trailing edge occupied by fans  \textless  \ 14.4 ($m$) & cont & 1 \\
 & Climb duration \textless  \ 1740 ($s $) & cont & 1 \\
 & Top of climb slope \textgreater \ 0.0108 ($rad$) & cont & 1 \\
 & \multicolumn{2}{l}{\textbf{Total number of constraints}} & {\textbf{5}} & \\
\hline
\end{tabular}
\label{tab:dragon}
\end{table}
For the position of the turboshaft layout, a categorical variable is defined with  2 associated levels (under the wing or behind).
To describe propulsive electric architectures, two main strategies can be considered.
\begin{description}
    \item[Categorical variable approach:] A single categorical variable is used to represent each possible combination of the number of generators (G) and the number of motors (M). This results in a total of 17 distinct choices (e.g., (2G,8M), (2G,12M), \ldots, (4G,40M),\ldots,(6G,36M)).
    \item[Hierarchical variable approach] where two hierarchical variables are used~\cite{halle2025distance}: a meta variable for the number of generators and a dependent variable for the number of motors, which is conditional on the number of generators selected.
\end{description}
The process minimizes fuel mass and leverages GP models adapted for mixed-categorical or hierarchical inputs within the SEGOMOE framework to handle the heterogeneous space efficiently~\cite{Saves2024HighDimMixed}. 
Different kernels are used: Gower Distance (GD), Continuous Relaxation (CR), and Homoscedastic Hypersphere (HH) for the first scenario, and Hierarchical (HIER) for the second scenario. These kernel-based GP optimization options are compared to the NSGA-II algorithm~\cite{nsga2} from the pymoo toolbox~\cite{pymoo} in Figure~\ref{fig:dragon}. The results demonstrate that BO with any of these kernel options is better suited to address such problems than evolutionary algorithms.

\subsection{Multi-fidelity \& multi-objective drone design} 

\begin{figure}[htp]
\begin{minipage}[b]{.6\linewidth}
\centering
\hspace{-1.25cm}
{
\includegraphics[width=7cm]
{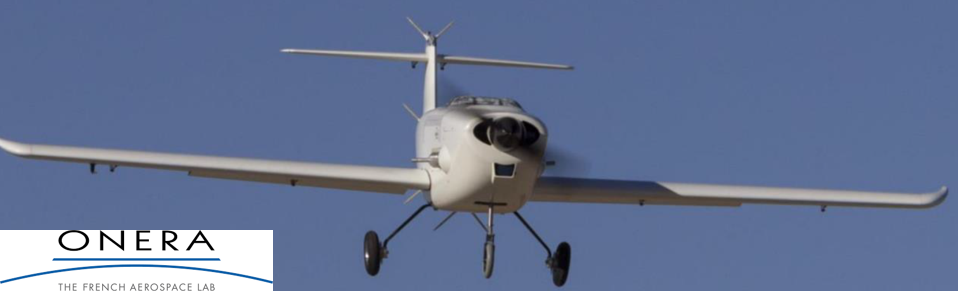}
}
\end{minipage}
\begin{minipage}[b]{.4\linewidth}
\centering 
\hspace{-1.5cm}
{\includegraphics[height=5cm,width=7cm]{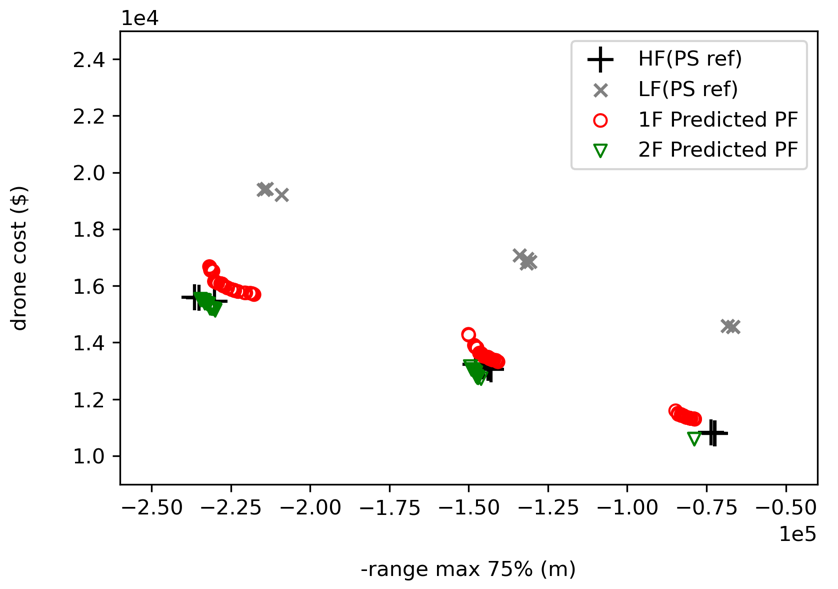}
}
\end{minipage}
\caption{Optimization results for the ``\texttt{K75}'' drone~\cite{SciTech_cat}: reference PF obtained with NSGA-II (on LF or HF models) compared to predicted PF with  BO.}
\label{fig:k75}
\end{figure}

We consider an electric version of a fixed-wing based drone: the Elecnor Deimos K75, represented in Figure~\ref{fig:k75}, available at ONERA. The analysis includes a suite of individual disciplinary models tailored to address
various aspects of fixed-wing drones, including aerodynamics, structure, autonomy, mission, cost
analysis, longitudinal stability, atmospheric conditions, and batteries. It is important to note
that the handling of these disciplines varied, including electric propulsion 
and mission complexity (simplified or relying on waypoints)~\cite{charayron2023advances}. Based on these variations, both a low-fidelity drone model and a high-fidelity one have been created.
Depending on the model fidelity, different tools are then encapsulated  within the OpenMDAO~\cite{gray2019openmdao} framework: OpenAeroStruct~\cite{jasa2018open}, DYMOS~\cite{Falck2021} for optimal control, financial cost~\cite{Charayron2023MultiFidelityMO,charayron2023advances}.

Table~\ref{tab:K75} presents the optimization problem with two objectives (the drone cost and the maximum range using only 75\% of the battery charge), two constraints (wing and tail failure), 19 continuous variables relative to the wing and tail geometry, 1 discrete variable for the number of battery blocks and 1 categorical variable for the 5 possible choices of wing material. The required quantities of interest are determined using the low-fidelity drone model and the high-fidelity drone
model characterized by a cost ratio of 262 in terms of CPU time.
A reference solution, obtained using the NSGA-II algorithm~\cite{Deb2002NSGAII,pymoo} with a large budget of 2000 high-fidelity (HF) evaluations (black crosses), is compared to the BO-predicted Pareto front in Figure~\ref{fig:k75}, which uses an equivalent budget of 20 HF evaluations (green triangles).
The Pareto front via BO obtained with only HF evaluations is represented by red circles, while the NSGA-II applied to low-fidelity (LF) evaluations is shown as gray crosses.
This drone design test case further demonstrates that the multi-fidelity methodology outperforms the mono-fidelity approach in terms of performance.
\begin{table}[htp]
\small
\centering
\vspace*{-0.3cm}
\caption{Definition of the ``\texttt{K75}'' optimization problem.}
\begin{tabular}{lllrp{2cm}}
& Function/variable & Nature & \# & Range\\
\hline
\hline
Minimize & Cost  & cont & 1 &\\
Maximize & Range & cont & 1 &\\
  & \multicolumn{2}{l}{\textbf{Total number of objectives}} & {\textbf{2}} & \\
\hline
w.r.t 
 & \multicolumn{2}{l}{Total  continuous variables} & 19 & \\
 \hline
 & Number of battery blocks & ordinal  & 3 & \{$2; 3; 4$\} \\
  &  Material & cat & 5 & \small\{fiberglass; Kevlar; carbon; HMcarbon;aluminium\} \\
 & \multicolumn{2}{l}{Total categorical variables} & 2 & \\
  \multicolumn{3}{l}{\textbf{Categorical strategy: total number of variables in relaxed dimension }} & {\textbf{27}} & \\
  \hline
subject to & $\text{wing}_{\text{failure}}(x) \leq 0$ & cont & 1 \\ &$\text{tail}_{\text{failure}}(x) \leq 0$ & cont & 1 \\
  & \multicolumn{2}{l}{\textbf{Total number of constraints}} & {\textbf{2}} & \\
\hline
\end{tabular}
\label{tab:K75}
\end{table}
\section{Conclusion}
Bayesian optimization, supported by flexible Gaussian process surrogate modeling and adaptive acquisition functions, provides an essential framework for solving highly constrained and complex black-box MDO problems. Recent developments successfully address key challenges: modeling of mixed and hierarchical variables~\cite{Saves2024SMT2, Saves2024HighDimMixed}, efficient solutions for multi-objective problems~\cite{Grapin2022ConstrainedMOBO, Charayron2023MultiFidelityMO}, handling complex non-linear constraints (UTB)~\cite{Priem2020UTB}, and managing hidden constraints (EFIFE)~\cite{Tfaily2024BayesianHC}. Furthermore, methods like EGORSE~\cite{Priem2024HighDimBO_AIAAJ}, CC-BO~\cite{Pretsch2025CooperativeComponentsBO}, and TREGO~\cite{diouane2022trego} extend BO applicability to high-dimensional problems, ensuring BO remains an efficient solution for designing innovative engineering systems. More complex problems, such as System of Systems (SoS), are currently under investigation.
\section*{Acknowledgments}
This work is part of the activities
of ONERA - ISAE - ENAC joint research group.
We would like to express our sincere gratitude to Rémy Priem and Rémy Charayron for their valuable contributions within SEGOMOE.  Special thanks to Luiz Tiberio Fernandez (ONERA), Jasper Bussemaker (DLR), Guiseppa Donelli (DLR), Lisa Pretsch (MTU), and Ali Tfaily (BOMBARDIER) for their insightful discussions and support throughout this work. 

Two open-source implementations of SEGOMOE core ideas for mono-objective and mono-fidelity optimizations are available with a full documentation including a large number of tutorials: the first one is a simple but effective implementation of the EGO BO algorithm within the SMT Python library\footnote{https://github.com/SMTorg/smt}~\cite{Saves2024SMT2}, the second one is an implementation of the constrained EGO BO algorithm within the EGObox library written in Rust\footnote{https://github.com/relf/egobox}~\cite{lafage2022egobox} (with Python bindings) aiming at better runtime performances and stability for the end-users.
\bibliography{main}
\end{document}